\documentclass[11pt]{article}
\usepackage{}

\usepackage[total={6in, 8.5in}]{geometry}
\usepackage{amsfonts}
\usepackage{amsthm}
\usepackage{amssymb}
\usepackage{amsmath}
\usepackage{enumerate}
\usepackage[
pdfauthor={ESYZ},
pdftitle={Toughness and spanning trees in K4mf graphs},
pdfstartview=XYZ,
bookmarks=true,
colorlinks=true,
linkcolor=blue,
urlcolor=blue,
citecolor=blue,
bookmarks=false,
linktocpage=true,
hyperindex=true
]{hyperref}

\makeatletter
\newcommand{\setword}[2]{%
	\phantomsection
	#1\def\@currentlabel{\unexpanded{#1}}\label{#2}%
}
\makeatother

\makeatletter
\renewcommand*\env@matrix[1][*\c@MaxMatrixCols c]{%
	\hskip -\arraycolsep
	\let\@ifnextchar\new@ifnextchar
	\array{#1}}
\makeatother

\usepackage{graphicx}
\usepackage[natural]{xcolor}
\usepackage{tikz}
\usepackage{tikz-3dplot}
\usetikzlibrary{shapes}
\usetikzlibrary{arrows,decorations.pathmorphing,backgrounds,positioning,fit,petri,automata}
\usetikzlibrary {positioning}
\usetikzlibrary{arrows}
\usepackage{tkz-euclide}

\usepackage{pgf,tikz,pgfplots}

\usepackage{mathrsfs}

\usetikzlibrary{arrows}

\usepackage{xcolor}
\long\def\ignore#1{}

\let\oldi\ignore

 \oldi{
\newtheorem{THM}{\textbf{Theorem}}[section]
\newtheorem{THMs}{\textbf{Theorem}}[section]
\newtheorem{DEF}[THM]{\textbf{Definition}}[section]
\newtheorem{LEM}[THM]{\textbf{Lemma}}
\newtheorem{CON}[THM]{\textbf{Conjecture}}
\newtheorem{PROP}[THM]{\textbf{Proposition}}
\newtheorem{COR}[THM]{\textbf{Corollary}}
\newtheorem{CORs}{\textbf{Corollary}}[section]
\newtheorem{PRO}[THM]{\textbf{Problem}}
\newcommand{\pf}{\textbf{Proof}.\quad}

\newtheorem{FAC}{\textbf{Fact}}
\newtheorem{REM}{\textbf{Remark}}
\newtheorem{OPR}{\textbf{Operation}}
\newtheorem{CLA}{\textbf{Claim}}[section]
\setcounter{CLA}{1}
 }

\newtheorem{THM}{Theorem}[section]

\newtheorem{LEM}[THM]{Lemma}
\newtheorem{CON}[THM]{Conjecture}

\newtheorem{COR}[THM]{Corollary}

\newtheorem{CLA}{Claim}[section]
\newcommand{\pf}{\textbf{Proof}.\quad}

\newtheorem*{LEM3}{\textbf{Lemma 3.3}}
\newtheorem*{LEM4}{\textbf{Lemma 2.8}}

\newcommand{\ve}{\varepsilon }

\newcommand{\D}{\Delta}
\newcommand{\de}{\delta}

\usepackage{thmtools}
\usepackage{thm-restate}

\newcommand{\arxiv}[1]{\href{http://arxiv.org/abs/#1}{\texttt{arXiv:#1}}}

\newcommand{\CC}{\mathcal{C}}

\newcommand{\pbar}{\overline{\varphi}}
 
\begin{document}
\title{The overfullness  of graphs with small minimum degree and large maximum degree}

\author{%
	Yan Cao\\
		Department of Mathematics, \\
	West Virginia University, Morgantown, WV 26506\\
	\texttt{yacao@mail.wvu.edu}
	\and
	Guantao Chen\thanks{This work was supported in part by NSF grant DMS-1855716.}\\
	Department of Mathematics and Statistics, \\
	Georgia State University, Atlanta, GA 30302\\
	\texttt{gchen@gsu.edu}%
	\and
	Guangming Jing\thanks{This work was supported in part by NSF grant DMS-2001130.}\\
	Department of Mathematics,\\
	Augusta University, Augusta, GA 30912\\
	\texttt{gjing@augusta.edu}%
	\and 
	Songling Shan\\
	Department of Mathematics, \\
	Illinois State Univeristy, Normal, IL 61790 \\
	\texttt{sshan12@ilstu.edu}
} 

\date{\today}
\maketitle

 \begin{abstract}
Given a simple graph $G$,    denote by $\Delta(G)$, $\delta(G)$, and $\chi'(G)$ the maximum degree, the minimum degree, and the chromatic index of $G$, respectively. 
We say $G$  is  \emph{$\Delta$-critical} if  $\chi'(G)=\Delta(G)+1$ and $\chi'(H)\le \Delta(G)$ for every proper subgraph $H$ of $G$; and  $G$ is \emph{overfull}  if $|E(G)|>\Delta \lfloor |V(G)|/2 \rfloor$. 
Since a maximum matching in $G$ can have size at most $\lfloor |V(G)|/2 \rfloor$, it follows that  $\chi'(G) = \Delta(G) +1$ if   $G$ is overfull.  
  Conversely, let $G$ be a $\Delta$-critical graph. 
 The well known overfull conjecture of Chetwynd and Hilton asserts that 
 $G$ is overfull provided $\Delta(G) > |V(G)|/3$.
  In this paper,  we show that any  $\Delta$-critical graph $G$ is overfull if   $\Delta(G)  - 7\delta(G)/4\ge  (3|V(G)|-17)/4$.


 \end{abstract}

\section{Introduction}
We will mainly use the notation from the book \cite{StiebSTF-Book}.  Graphs in this paper are simple. The
 vertex set and the  edge set of a graph $G$ are denoted by $V(G)$ and $E(G)$, respectively. We denote by $n(G)$, $\D(G)$ and $\de(G)$ the order, maximum degree and minimum degree of $G$, respectively.  When $G$ there is no risk of confusion, we simply use $n$, $\D$ and $\de$. 
For two integers $p,q$, let $[p,q]=\{i\in \mathbb{Z}:  p\le i \le q\}$.  
 
 An {\it edge $k$-coloring\/} or simply a \emph{$k$-coloring} of a graph $G$ is a mapping $\varphi$ from $E(G)$ to the set of integers
$[1,k]$, called {\it colors\/}, such that  no two adjacent edges receive the same color with respect to $\varphi$.  
The {\it chromatic index\/} of $G$, denoted $\chi'(G)$, is defined to be the smallest integer $k$ so that $G$ has an edge $k$-coloring.  
We denote by $\CC^k(G)$ the set of all edge $k$-colorings of $G$.  
A  graph $G$  is \emph{$\Delta$-critical} if  $\chi'(G)=\Delta+1$ and $\chi'(H)\le \Delta$ for every proper subgraph $H$ of $G$. 
In 1960's, Vizing~\cite{Vizing-2-classes} and, independently,  Gupta~\cite{Gupta-67} proved that $\Delta \le \chi'(G) \le \Delta+1$. 
This leads to a natural classification of graphs. Following Fiorini and Wilson~\cite{fw},  we say a graph $G$ is of {\it class 1} if $\chi'(G) = \Delta$ and of \emph{class 2} if $\chi'(G) = \Delta+1$.  Holyer~\cite{Holyer} showed that it is NP-complete to determine whether an arbitrary graph is of class 1.  
Nevertheless, if a graph $G$ has too many edges in comparing with its maximum degree, i.e., $|E(G)|>\Delta \lfloor n/2\rfloor$,  then we have to color $E(G)$ using exactly  $\Delta+1$ colors. Such graphs are  \emph{overfull}. Although the converse is not true in general, finding sufficient conditions for a class 2 graph to be overfull has been the major focus in graph edge chromatic theory.

Applying Edmonds' matching polytope theorem, Seymour~\cite{seymour79}  showed  that whether a graph  $G$ contains an overfull subgraph of maximum degree $\Delta(G)$ can be determined in polynomial time.  A number of long-standing conjectures listed in {\it Twenty Pretty Edge Coloring Conjectures} in~\cite{StiebSTF-Book} lie in deciding when a $\Delta$-critical graph is overfull.   Chetwynd and  Hilton~\cite{MR848854,MR975994},  in 1986, proposed the following 
conjecture. 
\begin{CON}[Overfull conjecture] \label{overfull-con}
If $G$  is a  class 2 graph  with  $\Delta(G)>\frac{n}{3}$,  then $G$ is contains an overfull subgraph $H$  with $\Delta(H)=\Delta(G)$. 
\end{CON}

 The degree condition that $\Delta>\frac{n}{3}$ in the conjecture above is best possible, as seen by the 
graph  obtained from the Petersen graph by deleting one vertex. 
If the overfull conjecture is true, then the NP-complete problem of 
determining the chromatic index becomes  polynomial-time solvable 
for graphs $G$ with $\Delta>\frac{n}{3}$. Furthermore, the overfull 
conjecture implies several other conjectures in edge colorings 
such as the 1-factorization conjecture. 
Despite its importance, very little is known about its truth.
It was confirmed only for graphs with $\Delta \ge  n-3$ by 
Chetwynd and   Hilton~\cite{MR975994} in 1989. By restricting the minimum degree, 
Plantholt~\cite{MR2082738} in 2004 showed that the overfull conjecture is affirmative  for 
graphs $G$ with  even order $n$ and minimum degree $\delta \ge \sqrt{7}n/3\approx 0.8819 n$. 
The 1-factorization conjecture is a special case of the overfull conjecture (when $G$ is regular with $\Delta\ge n/2$), which in 2013 
was confirmed for large $n$ by Csaba, K\"uhn, Lo, Osthus, Treglown
~\cite{MR3545109}. The overfull conjecture is still wide open in general, and it seems
extremely difficult even for graphs $G$ with  $\Delta \ge n-4$. 

Since every class 2 graph of maximum degree $\Delta$ contains a $\Delta$-critical subgraph, the overfull conjecture is equivalent to 
saying that every $\Delta$-critical graph with  $\Delta(G)>\frac{n}{3}$ is overfull. 
In this paper, we prove the following result 
 and hope to apply it 
to verify the overfull conjecture for graphs  with $\Delta  \ge (1-\varepsilon)n$ for a small positive real number $\varepsilon >0$.

\begin{THM}\label{thm:overfull-min-degree}
Let $G$ be a  $\Delta$-critical graph of order $n$.  If 
  $\Delta(G) -\frac{7\delta(G)}{4} \ge \frac{3n-17}4$,  then $G$ is overfull. 
\end{THM}

\begin{COR}
Let $0<\ve<1/7$ be any given real number and $G$ be a $\Delta$-critical graph of order $n$. 
Then  $G$ is overfull provided that 
$\delta(G) \le  \ve n $ and 
 $\Delta \ge \frac{3n-17+7\ve n}{4}$. 
\end{COR}

The remainder of this paper is organized as follows. In next section, we introduce some preliminaries on edge colorings such as Vizing fans and Kierstead paths. In Section 3, we prove Theorem~\ref{thm:overfull-min-degree}.
In the last section, we prove one of the main lemmas. 


\section{Preliminaries}\label{lemma}

In this section, we introduce the fundamental tools such as multifans and Kierstead paths 
in edge colorings, and list some new developments built on these tools. 
Let $G$ be a graph and  $e\in E(G)$.   Denote by $G-e$
the graph obtained from $G$ by deleting the edge $e$. 
  A vertex is called a \emph{$k$-vertex}  if its degree is $k$;   and a neighbor of a vertex $v$ is a 
  \emph{$k$-neighbor} if it  is a $k$-vertex in $G$. 
  
An edge $e\in E(G)$ is a \emph{critical edge} of $G$ if $\chi'(G-e)<\chi'(G)$. 
It is not hard to see that  a connected graph
  is $\D$-critical if and only if  every of its edge is critical.  
Critical graphs are useful since they have much more graph structural properties  than general class 2 graphs do. For 
example,  Vizing's Adjacency Lemma (VAL) from 1965~\cite{Vizing-2-classes} is a useful tool in counting the number of $\D$-neighbors of a vertex. 
\begin{LEM}[VAL]Let $G$ be a class 2 graph. If $xy$ is a critical edge of $G$, then $x$ has at least $\Delta-d_G(y)+1$ $\Delta$-neighbors from $V(G)\setminus \{y\}$.
	\label{thm:val}
\end{LEM}

Let $G$ be a graph, $e\in E(G)$ and 
$\varphi\in \CC^k(G-e)$ for  some integer $k\ge 0$. 
For any $v\in V(G)$, the set of colors \emph{present} at $v$ is 
$\varphi(v)=\{\varphi(f)\,:\, \text{$f$ is incident to $v$}\}$, and the set of colors \emph{missing} at $v$ is $\pbar(v)=[1,k]\setminus\varphi(v)$.  If $|\pbar(v)|=1$, we will also use $\pbar(v)$ to denote the  color that is missing at $x$.
For a vertex set $X\subseteq V(G)$,  define 
$
\pbar(X)=\bigcup _{v\in X} \pbar(v).
$ The set $X$ is called \emph{elementary} with respect to $\varphi$  or simply \emph{$\varphi$-elementary} if $\pbar(u)\cap \pbar(v)=\emptyset$
for every two distinct vertices $u,v\in X$.   Moreover,  we sometimes just say that $X$ 
is elementary if the  edge coloring is understood. For two distinct colors $\alpha,\beta \in [1,\Delta]$, the components of the subgraph induced by edges with colors $\alpha$ or $\beta$ are called $(\alpha, \beta)$-chains. Clearly, each $(\alpha, \beta)$-chain is either a path or an even cycle.  For an $(\alpha,\beta)$-chain  $P$, if it is a path with an endvertex $x$, 
we also denote it by $P_x(\alpha,\beta,\varphi)$ to stress the endvertex $x$. 
  If we interchange the colors $\alpha$ and $\beta$
 on an $(\alpha,\beta)$-chain $C$ of $G$, we get a new edge $k$-coloring  of $G$,  which is denoted by  $\varphi/C$.  This operation is a \emph{Kempe change}. 
 If $C$ is a path with one endvertex $x$, we may also refer this Kempe change as 
 ``do an $(\alpha,\beta)$-swap at $x$.'' 
 
 Let  $x,y\in V(G)$.   If $x$ and $y$
 are contained in a same  $(\alpha,\beta)$-chain of $G$ with respect to $\varphi$, we say $x$ 
 and $y$ are \emph{$(\alpha,\beta)$-linked} with respect to $\varphi$.
 Otherwise, $x$ and $y$ are \emph{$(\alpha,\beta)$-unlinked} with respect to $\varphi$. Without specifying $\varphi$, when we just say  $x$ and $y$ are $(\alpha,\beta)$-linked or $x$ and $y$ are $(\alpha,\beta)$-unlinked, we mean they are linked or unlinked with respect to the current edge coloring.

A \emph{multi-fan} at $x$ with respect to edge $e=xy \in E(G)$ and coloring $\varphi\in \mathcal{C}^\Delta(G-e)$ is a sequence $F = (x, e_1,y_1,\ldots, e_p,y_p)$ with $p \ge 1$ consisting of edges $e_1,e_2,\ldots,e_p$ and distinct vertices $x,y_1,y_2,\ldots,y_p$ satisfying the following two conditions:
\begin{itemize}
\item  $e_1=e$ and $e_i =  xy_i$ for $i=1, \ldots, p$,
\item for every edge $e_i$ with $i\in[2,p]$, there is  
$j\in [1,i-1]$ such that $\varphi(e_i)\in \pbar(y_j)$.
\end{itemize}
Denote by $V(F)$ the set of vertices contained in $F$. 

Notice that multi-fan is slightly more general than Vizing-fan which requires $j=i-1$ in the second condition. The following lemma shows that the vertex set of a multi-fan is elementary. The proof can be found in the book~\cite{StiebSTF-Book}.

\begin{LEM}{\em [Stiebitz, Scheide, Toft and Favrholdt~\cite{StiebSTF-Book}]}\label{vf}
 Let $G$ be a class 2 graph with maximum degree $\Delta$, $e_1\in E(G)$, $\varphi\in \CC^\Delta(G-e_1)$, and let  $F = (x, e_1, y_1, \ldots, e_p, y_p)$ be a multi-fan at $x$ with respect to $e_1$ and $\varphi$. Then the following statements hold:

(a) $\{x,y_1,y_2,\ldots,y_p\}$ is $\varphi$-elementary.

(b) For any $\alpha \in \pbar(x)$ and $\beta \in \pbar(y_i)$ for  some $i \in[1,p]$, it holds that  $P_x(\alpha,\beta,\varphi)=P_{y_i}(\alpha,\beta,\varphi)$.

\end{LEM}

Let a multifan $F = (x, e_1, y_1, \ldots, e_p, y_p)$ be at $x$ with respect to $e_1$ and $\varphi$, and $y_{\ell_1}, \ldots, y_{\ell _k}$
be a subsequence of $y_1, \ldots, y_p$.  We call $y_{\ell_1},y_{\ell_2}, \ldots, y_{\ell_k}$ an  \emph{$\alpha$}-sequence with respect to $\varphi$ and $F$ if the following holds:
$$
\varphi(xy_{\ell_1})= \alpha\in \pbar(y_1),  \quad \varphi(xy_{\ell_i})\in \pbar(y_{\ell_{i-1}}), \quad  i\in [2,k].
$$
A vertex in an $\alpha$-sequence is called an \emph{$\alpha$-inducing vertex} with respect to $\varphi$ and $F$, and a missing color at an $\alpha$-inducing vertex is called an \emph{$\alpha$-inducing color}. For convenience, $\alpha$ itself is also an $\alpha$-inducing color. We say a color $\beta$ is {\it induced by} $\alpha$ if $\beta$ is $\alpha$-inducing. By Lemma~\ref{vf} (a) and the definition of multifan, each color in $\pbar(V(F))$ is induced by a unique color in $\pbar(y_1)$. Also if $\alpha_1,\alpha_2$ are two distinct colors in $\pbar(y_1)$, then an $\alpha_1$-sequence is disjoint with an $\alpha_2$-sequence. For two distinct $\alpha$-inducing colors $\beta$ and $\delta$, we write {$\mathit \delta \prec \beta$} if there exists an $\alpha$-sequence $y_{\ell_1},y_{\ell_2}, \ldots, y_{\ell_k}$ such that $\delta\in\pbar(y_{\ell_i})$, $\beta\in\pbar(s_{\ell_j})$ and $i<j$. For convenience, $\alpha\prec\beta$ for any $\alpha$-inducing color $\beta\not=\alpha$.
As a consequence of Lemma~\ref{vf} (a), we have the following properties for a multifan. 
A proof of the result can be found in~\cite[Lemma 3.2]{HZ}. 
\begin{LEM}
	\label{thm:vizing-fan2}
	Let $G$ be a Class 2 graph and $F = (x, e_1, y_1, \ldots, e_p, y_p)$  be a multifan with respect to a critical edge $e_1=xy_1$ and a coloring $\varphi\in \CC^\Delta(G-e_1)$. For two colors $\delta\in \pbar(y_i)$ and $\lambda\in \pbar(y_j)$ with  $i,j\in [1,p]$ and $i\ne j$, the following statements  hold.
	\begin{enumerate}[(a)]
		\item If $\delta$ and $\lambda$ are induced by different colors, then $y_i$ and $y_j$ are $(\delta, \lambda)$-linked with respect to $\varphi$. 
		\label{thm:vizing-fan2-a}
		\item If $\delta$ and $\lambda$ are induced by the same color, $\delta\prec\lambda$, and $y_i$ and $y_j$ are $(\delta, \lambda)$-unlinked with respect to $\varphi$, 
		then $x\in P_{y_j}(\lambda, \delta, \varphi)$.  	\label{thm:vizing-fan2-b}
	\end{enumerate}
	
\end{LEM}

Let $G$ be a graph, $e=v_0v_1\in E(G)$, and  $\varphi\in \CC^k(G-e)$ for some integer $k\ge 0$.
A \emph{Kierstead path}  with respect to $e$ and $\varphi$
is a sequence $K=(v_0, v_0v_1, v_1, v_1v_2, v_2, \ldots, v_{p-1}, v_{p-1}v_p,  v_p)$ with $p\geq 1$ consisting of  distinct vertices $v_0,v_1, \ldots , v_p$ and distinct edges $v_0v_1, v_1v_2,\ldots, v_{p-1}v_p$ satisfying the following condition:
\begin{enumerate}[(K1)]
	\item For every edge $v_{i-1}v_i$ with $i\in [2,p]$,  there exists $j\in [0,i-2]$ such that 
	$\varphi(v_{i-1}v_i)\in \pbar(v_j)$. 
\end{enumerate}

Clearly, a Kierstead path with at most 3 vertices is a multifan, and so the set of its vertices is elementary. The following  result shows that the set of the vertices of Kierstead paths with $4$ vertices is ``nearly'' elementary.  (See Theorem 3.3 from~\cite{StiebSTF-Book}.) 

\begin{LEM}[]\label{Lemma:kierstead path1}
	Let $G$ be a Class 2 graph,
	$e=v_0v_1\in E(G)$ be a critical edge, and $\varphi\in \CC^\Delta(G-e)$. If $K=(v_0, v_0v_1, v_1, v_1v_2,  v_2, v_2v_3, v_3)$ is a Kierstead path with respect to $e$
	and $\varphi$, then the following statements hold.
	\begin{enumerate}[(a)]
		\item If $\min\{d_G(v_1), d_G(v_2)\}<\Delta$, then $V(K)$ is $\varphi$-elementary.
		\item  $|\pbar(v_3)\cap (\pbar(v_0)\cup \pbar(v_1))|\le 1$. 
	\end{enumerate}

\end{LEM}

\begin{LEM}\label{degreek}
	Let   $G$ be an $n$-vertex  $\Delta$-critical graph  and let $a\in V(G)$.   If  $d(a) \le \frac{2\D -n +2}3$,   then for each $v\in V(G)\setminus \{a\}$,
	either $d(v) \ge \Delta-d(a)+1$ or $d(v) \le n-\Delta+2d(a)-6$. Furthermore,  if $d(v) \ge \Delta-d(a)+1$, 
	then for any $b\in N(a)$ with $d(b)=\Delta$ and $\varphi\in \CC^\Delta(G-ab)$,   $|\pbar(v)\cap (\pbar(a)\cup \pbar(b))| \le 1$.
\end{LEM}

\pf Let $k =d(a)$. Assume to the contrary that there exists $v\in V(G)\setminus \{a\}$
such that $n-\Delta+2k-5 \le d(v) \le \Delta-k$. 
Let $b\in N(a)$ with $d(b)=\Delta$ and $\varphi\in \CC^\Delta(G-ab)$.  
As $d(b) = \Delta$, we have $| N(b)\setminus \{a\}| = \Delta-1$.
Since $d(v) \le \Delta-k$, by VAL, $a\not\in N(v)$. 
As $b,v\not\in N(b)\cap N(v)$ and $d(v)  \ge  n-\Delta+2k-5$, it follows that 
\begin{eqnarray*}
|N(v)\cap (N(b)\setminus\{a\})| &=&|((N(v) \cap (N(b)\setminus \{a\}))\setminus \{b,v\}|\\
&=&  |N(v)|+|N(b)\setminus \{a\}|-(n-3)-|N(b)\cap \{v\}|-|N(v)\cap \{b\}| \\
&\ge & 2k-3-|N(b)\cap \{v\}|-|N(v)\cap \{b\}|.
\end{eqnarray*}
Note that $|[1,\Delta]\setminus (\pbar(a)\cup \pbar(b))|=k-2 $. If $bv\not\in E(G)$, we can find 
a vertex $u\in N(v)\cap \big(N(b)\setminus \{a\}\big)$ 
such that $\varphi(bu), \varphi(vu) \in \pbar(a)\cup \pbar(b)$. 
If $bv\in E(G)$, then as $d(v) \le \Delta-k$, by Lemma~\ref{vf} (a), 
we must have $\varphi(bv) \not\in \pbar(a)\cup \pbar(b)$. 
Thus there are at most $2(k-3)$ edges between $\{b,c\}$ and $N(v)\cap (N(b)\setminus\{a\})$
colored by colors from $[1,\Delta]\setminus (\pbar(a)\cup \pbar(b))$. 
Thus again we can find 
a vertex $u\in N(v)\cap \big(N(b)\setminus \{a\}\big)$ 
such that $\varphi(bu), \varphi(vu) \in \pbar(a)\cup \pbar(b)$. 
Therefore, $K=(a,ab,b,bu,u,uv,v)$ is a Kierstead path. 
Since $d(v) \le \Delta-k$, $|\pbar(v)| \ge k$.
This implies that 
$$
|\pbar(v)\cap (\pbar(a)\cup \pbar(b)) | \ge 2, 
$$
showing a contradiction to Lemma~\ref{Lemma:kierstead path1} (b). 

For the second part, since $| N(b)\setminus \{a\}| =\Delta-1$ and $d(v) \ge \Delta-k+1$, 
$$
|N(v)\cap \big(N(b)\setminus \{a\}\big)|  \ge 2\Delta -n-k \ge  2k-2,
$$
as $k=d(a) \le \frac{2\D -n +2}3$. 
The rest of the proof then follows  the same argument as for the first part, and so we omit it. 
\qed 

Lemmas~\ref{lem:fork} to \ref{lem:kite} below aim to 
study the ``elementary properties'' of small tree-like structures 
that contain Kierstead paths as substructures. 

Let $G$ be a $\Delta$-critical graph, $ab\in E(G)$, and $\varphi \in \CC^\Delta(G-ab)$. 
A {\it fork}  $H$  with respect to $\varphi$ is a graph with 
$$V(H)=\{a,b,u,s_1,s_2,t_1,t_2\} \quad \text{and}\quad E(H)=\{ab,bu,us_1,us_2, s_1t_1,s_2t_2\}$$ 
such that $\varphi(bu)\in \pbar(a)$, $\varphi(us_1), \varphi(us_2) \in \pbar(a)\cup \pbar(b)$,  and $\varphi(s_1t_1)\in (\pbar(a)\cup \pbar(b))\cap \pbar(t_2) $
and $\varphi(s_2t_2)\in (\pbar(a)\cup \pbar(b))\cap \pbar(t_1)$. 
Fork  was defined in~\cite{av2} and  
it was shown in~\cite[Proposition B]{av2} that a fork can not exist in a $\Delta$-critical graph if 
the degree sum of $a$, $t_1$ and $t_2$ is small.

\begin{LEM}\label{lem:fork}
	Let $G$ be a $\Delta$-critical graph, $ab\in E(G)$, and $\{u,s_1,s_2, t_1,t_2\}\subseteq V(G)$. 
	If $\Delta\ge d_G(a)+d_G(t_1)+d_G(t_2)+1$, 
	then for any $\varphi\in \CC^\Delta(G-ab)$, $G$ does not contain a fork on $\{a,b,u,s_1,s_2,t_1,t_2\}$ with respect to $\varphi$.   
\end{LEM}

A {\it short-kite}  $H$ is a graph with 
$$V(H)=\{a,b,c,u,x,y\} \quad \text{and}\quad E(H)=\{ab,ac,bu,cu,ux,uy\}.$$ 
The lemma below was proven in~\cite{2103.05171}.  

\begin{LEM}\label{lemma:class2-with-fullDpair2}
	Let  $G$ be  a class 2 graph, 
	$H\subseteq G$ 
	be a short-kite with $V(H)=\{a,b,c,u,x,y\}$, and let $\varphi\in \CC^\Delta(G-ab)$. 
	Suppose $$K=(a,ab,b,bu,u,ux,x) \quad \text{and} \quad K^*=(b,ab,a,ac,c,cu,u,uy,y)$$
	are two Kierstead path with respect to $ab$ and $\varphi$.  
	If $\pbar(x)\cap ( \pbar(a)\cup \pbar(b))\ne \emptyset$ and $\pbar(y)\cap ( \pbar(a)\cup \pbar(b))\ne \emptyset$,  then $\max\{d(x),d(y)\}=\Delta $. 
\end{LEM}

A {\it kite}  $H$ is a graph with 
$$V(H)=\{a,b,c,u,s_1,s_2,t_1,t_2\} \quad \text{and}\quad E(H)=\{ab,ac,bu,cu,us_1,us_2, s_1t_1,s_2t_2\}.$$ 
\begin{LEM}\label{lem:kite}
	Let $G$ be a class 2 graph, $H\subseteq G$ 
	be a kite with $V(H)=\{a,b,c,u,s_1,s_2,t_1,t_2\}$, and let $\varphi\in \CC^\Delta(G-ab)$. 
	Suppose $$K=(a,ab,b,bu,u,us_1, s_1,s_1t_1,t_1) \quad \text{and} \quad K^*=(b,ab,a,ac,c,cu,u,us_2, s_2,s_2t_2,t_2)$$
	are two Kierstead paths with respect to $ab$ and $\varphi$. 
	If $\varphi(s_1t_1)=\varphi(s_2t_2)$, 
	then $|\pbar(t_1)\cap \pbar(t_2) \cap ( \pbar(a)\cup \pbar(b))|\le 4$.  
\end{LEM}
The proof of Lemma~\ref{lem:kite} will be given in the last section.

\section{Proof of Theorem~\ref{thm:overfull-min-degree}}\label{lemma}

  Since all vertices not missing  a given color $\alpha$
  are saturated by the matching that consists of all edges colored by $\alpha$ in $G$, we have the Parity Lemma below, which has appeared in many papers, for example, see~\cite[Lemma 2.1]{MR2028248}.  
  \begin{LEM}[Parity Lemma]
  	Let $G$ be an $n$-vertex multigraph and $\varphi\in \CC^k(G)$ for some integer $k\ge \Delta(G)$. 
  	Then for any color $\alpha\in [1,\Delta]$, 
  	$|\{v\in V(G): \alpha\in \pbar(v)\}| \equiv n \pmod{2}$. 
  \end{LEM}

\proof[Proof of Theorem~\ref{thm:overfull-min-degree}]
Since the overfull conjecture is true for $\Delta \ge  n-3$ by 
Chetwynd and   Hilton~\cite{MR975994}, we may assume that $\Delta \le n-4$.
Thus $\Delta(G) -\frac{7\delta(G)}{4} \ge \frac{3n-17}4$ implies that 
$\delta(G) \le \frac{n+1}{7}$. As $\delta(G) \ge 2$, it follows that $n\ge 13$ and $\Delta(G) \ge \frac{3n-3}{4}$.  

Choose $a\in V(G)$ such that $d(a)=\delta(G)=:k$. As $\Delta \ge \frac{3n-3}{4}$, 
$ \frac{2\Delta-n+2}{3} \ge \frac{n+1}{6}>k$. 
By Lemma~\ref{degreek}, for every 
$v\in V(G)\setminus\{a\}$, we have either 
$$
d(v) \ge \Delta-k+1 \quad \text{or} \quad d(v) \le (n-\Delta)+2k-6. 
$$
Let $b\in N(a)$ with $d(b)=\Delta$ and $\varphi\in \CC^\Delta(G-ab)$. 
We consider two cases below to finish the proof. 

\medskip 

{\bf \noindent Case 1}: there exist two distinct vertices $t_1,t_2 \in V(G)\setminus \{a\}$ such that $d(t_1), d(t_2)\le (n-\Delta)+2k-6$.

\medskip 

Let $s_1\in N(t_1)$ and $s_2\in N(t_2)$ be any two distinct vertices such that $\varphi(s_1t_1), \varphi(s_2t_2) \in \pbar(a)\cup \pbar(b)$. Such vertices $s_1, s_2$ exists since $d(a)=k=\delta(G)$ and $|[1,\Delta]\setminus (\pbar(a)\cup \pbar(b))|=k-2$,
there are at least two edges incident to each of $t_1$ and $t_2$ colored by colors from $\pbar(a)\cup \pbar(b)$. 
Since $G$ is $\Delta$-critical, VAL implies that for 
every $s_i\in N(t_i)$, $d(s_i)\ge \Delta+2-(n-\Delta+2k-6)>n-\Delta+2k-6$, 
where the last inequality follows as $\Delta -\frac{7k}{4} \ge \frac{3n-17}4$
and $n\ge 13$.
 
By Lemma~\ref{degreek}, $d(s_i)\ge \Delta-k+1$. In fact we can assume that $d(s_i)\ge \Delta-k+2$. Since if $d(s_i)=\Delta-k+1$, then $t_i$ has at least $k-1$ neighbors of degree $\Delta$ by VAL. Note that $k-1>k-2=|[1,\Delta]\backslash (\pbar(a)\cup \pbar(b))|$, we may choose $s_i$ such that $\varphi(s_it_i)\in \pbar(a)\cup \pbar(b)$ with $d(s_i)=\Delta\ge \Delta-k+2$. Thus we assume that $d(s_i)\ge \Delta-k+2$. Let $c\in N(a)$ such that $\varphi(ac) \in \pbar(b)$. By VAL, $d(c) \ge \Delta-k+2$. 
Thus, since $b,c\not\in N(b)\cap N(c)$ and $s_1,s_2\not\in N(s_1)\cap N(s_2)$, 
\begin{eqnarray}\label{common-neg}
&&|N(b)\cap N(c)\cap N(s_1)\cap N(s_2)|\\ \nonumber 
&=& |(N(b)\cap N(c)\cap N(s_1)\cap N(s_2))\backslash\{b,c,s_1,s_2\}| \nonumber \\
&\ge & |N(s_1)\cap N(s_2)|-|N(s_1)\cap N(s_2)\cap \{b,c\}| \nonumber \\ && +|N(b)\cap N(c)|-|N(b)\cap N(c)\cap \{s_1,s_2\}|-(n-4) \nonumber \\
&\ge & (2\Delta-n-2k+4)-|N(s_1)\cap N(s_2)\cap \{b,c\}| \nonumber \\
&& +(2\Delta-n-k+2)-|N(b)\cap N(c)\cap \{s_1,s_2\}|-n+4 \nonumber \\
&=& 4\Delta-3n-3k+10-|N(s_1)\cap N(s_2)\cap \{b,c\}|\nonumber-|N(b)\cap N(c)\cap \{s_1,s_2\}| \nonumber \\
&\ge& (k-2  -|N(s_1)\cap\{b,c\}|)+(k-2-|N(s_2)\cap\{b,c\}|)+\nonumber \\ &&(k-2-|N(b)\cap \{s_1,s_2\}|)+(k-2-|N(c)\cap\{s_1,s_2\}|)+1, \nonumber 
\end{eqnarray}
as $\Delta\ge 3n/4+(7k-17)/4$. Note that $|(\pbar(a)\cup \pbar(b)) \cap \pbar(t_i)| \ge 2$ under any $\Delta$-coloring of $G-ab$ since $d(t_i) \le (n-\Delta)+2k-6 \le \Delta-k$. Therefore if $s_ib\in E(G)$, then $\varphi(s_ib)\notin \pbar(a)\cup \pbar(b)$ by Lemma~\ref{Lemma:kierstead path1} (b). By considering the coloring $\varphi'$ obtained from $\varphi$ by coloring $ab$ with $\varphi(ac)$ and uncoloring $ac$, we also have that if $s_ic\in E(G)$, then $\varphi(s_ic)\notin \pbar(a)\cup \pbar(b)$.
Since $|[1,\Delta]\setminus (\pbar(a)\cup \pbar(b))|=k-2 $, {\bf (1)} implies that there exists  $u\in N(b)\cap N(c)\cap N(s_1)\cap N(s_2)$
such that $\varphi(ub), \varphi(uc), \varphi(us_1),\varphi(us_2) \in \pbar(a)\cup \pbar(b)$.  As both $\{a,b,t_i\}$  and $\{b,a, s_i, t_i\}$
are not elementary,  it follows that $u\notin \{a,t_1,t_2\}$
by Lemma~\ref{vf} and Lemma~\ref{Lemma:kierstead path1}.
Thus  $H$ with $V(H)=\{a,b,c,u,s_1,s_2,t_1,t_2\}$
is a kite, and both $$K=(a,ab,b,bu,u,us_1, s_1,s_1t_1,t_1) \quad \text{and} \quad K^*=(b,ab,a,ac,c,cu,u,us_2, s_2,s_2t_2,t_2)$$
are  Kierstead paths with respect to $ab$ and $\varphi$.

If there exist $s_1\in N(t_1)$ and $s_2\in N(t_2)$ such that $\varphi(s_1t_1), \varphi(s_2t_2) \in \pbar(a)\cup \pbar(b)$ and $\varphi(s_1t_1)=\varphi(s_2t_2)$,
then let  $\Gamma= \pbar(t_1)\cap \pbar(t_2)$.  Since $|[1,\Delta]\setminus (\pbar(a)\cup \pbar(b))|=k-2$, $|\Gamma \cap (\pbar(a)\cup \pbar(b))|\ge |\Gamma|-k+2$. 
By the assumption of this case, 
$|\Gamma|\ge \Delta-2(n-\Delta+2k-6)=3\Delta-2n-4k+12$. As  
$\Delta \ge \frac{3n+7k-17}{4}$, $n\ge 13$ and $k\ge 2$,   we get $$|\Gamma \cap (\pbar(a)\cup \pbar(b))|\ge |\Gamma|-k+2 \ge 3\Delta-2n-5k+14 \ge \frac{n+k+5}{4}\ge 5>4,$$
contradicting  Lemma~\ref{lem:kite}.

Thus for any $s_1\in N(t_1)$ and $s_2\in N(t_2)$ such that $\varphi(s_1t_1), \varphi(s_2t_2) \in \pbar(a)\cup \pbar(b)$, it follows that $\varphi(s_1t_1) \ne \varphi(s_2t_2)$. Thus $\varphi(s_1t_1)\in \pbar(t_2)$ and $\varphi(s_2t_2)\in \pbar(t_1)$. 
Hence $H^*$ with $V(H^*)=\{a,b,u,s_1,s_2,t_1,t_2\}$
is a fork.
However,   $d(a)+d(t_1)+d(t_2)\le k+2(n-\Delta+2k-6)=2n-2\Delta+5k-12 <\Delta$, as $n\ge 13$
and $\Delta\ge 3n/4+(7k-17)/4$, contradicting Lemma~\ref{lem:fork}.

\medskip 

{\bf \noindent Case 2}: there exists at most one vertex $t \in V(G)\setminus \{a\}$ such that $d(t)\le (n-\Delta)+2k-6$.

\medskip

Assume that $G$ is not overfull. We show that we can always find two distinct vertices $x,y\in V(G)\setminus \{a\}$
with 
$$ \Delta-k+1\le d(x), d(y)<\Delta, \quad  \pbar(x)\cap (\pbar(a)\cup \pbar(b))\ne \emptyset, \quad \text{and}\quad 
\pbar(y)\cap (\pbar(a)\cup \pbar(b))\ne \emptyset. \quad (*)$$ Then we will apply Lemma~\ref{lemma:class2-with-fullDpair2} to reach a contradiction. To see this, let $c\in N(a)$ such that $\varphi(ac) \in \pbar(b)$. By VAL, $d(c) \ge \Delta-k+2$.
As $d(x), d(y) \ge \Delta-k+1$, by  exactly the same argument as in ~\eqref{common-neg}, we can find $u\in N(b)\cap N(c)\cap N(x)\cap N(y)$
such that $u\not\in \{a,x,y\}$ and $\varphi(ub), \varphi(uc), \varphi(ux),\varphi(uy) \in \pbar(a)\cup \pbar(b)$. 
Then  $H$ with $V(H)=\{a,b,c,u,x,y\}$
is a short-kite, and both $$K=(a,ab,b,bu,u,ux, x) \quad \text{and} \quad K^*=(b,ab,a,ac,c,cu,u,uy, y)$$
are  Kierstead paths with respect to $ab$ and $\varphi$. As $\max\{d(x), d(y)\}<\Delta$, we achieve a contradiction to Lemma~\ref{lemma:class2-with-fullDpair2}. Thus we only show the existence of the vertices $x$ and $y$ below. 

Assume first that $n$ is odd. As $G$ is not overfull, some color in $[1,\Delta]$ is missed at at least three distinct vertices 
of $G$ by the Parity Lemma. 
If there exists exactly one vertex $t \in V(G)\setminus \{a\}$ such that $d(t)\le (n-\Delta)+2k-6$,
then each color in $\pbar(t)\cap (\pbar(a)\cup \pbar(b))$ is missed at another vertex from $V(G)\setminus \{a,b,t\}$. 
As $|\pbar(t)\cap (\pbar(a)\cup \pbar(b))|\ge 2\Delta-n-2k+6-(k-2)=2\Delta-n-3k+8>k$,   and $d(v)\ge \Delta-k+1$
for every $v\in V(G)\setminus \{a,b,t\}$, we can find $x,y\in V(G)\setminus \{a\}$ 
with the desired property as in $(*)$. 
Thus for every vertex $v\in V(G)\setminus \{a,b\}$, $d(v)\ge \Delta-k+1$. 
Let $\delta\in [1,\Delta]$ such that $\delta$ is missing at at least three vertices, say $u,v,w$, from $V(G)$. 
If $\delta \in \pbar(a)\cup \pbar(b)$, as $\pbar(a)\cup \pbar(b)$ is $\varphi$-elementary,  then letting $x,y\in \{u,v,w\}\setminus \{a,b\}$ will give us a desired choice. 
So we assume $\delta \not\in \pbar(a)\cup \pbar(b)$. 
Let $\alpha\in \pbar(a)\cup \pbar(b)$, say, w.l.g., that  $\alpha\in \pbar(a)$. 
As at most one of $u,v,w$, say $w$ is $(\alpha,\delta)$-linked with $a$, we  do an $(\alpha,\delta)$-swap at the other 
two vertices $u$ and $v$. Call the new coloring still $\varphi$, now $u$ and $v$ can play the role of $x$
and $y$.  

Assume now  that $n$ is even. As $G$ is not overfull,   each color in $[1,\Delta]$ is missing at an even number of vertices 
of $G$ by the Parity Lemma.  In particular, we have fact (a): each color in $\pbar(a)\cup \pbar(b)$
is missing at a vertex from $V(G)\setminus \{a,b\}$.  By the second part of Lemma~\ref{degreek}, we have fact (b): for every $v\in V(G)\setminus \{a,b\}$, if $  \Delta-k+1 \le d(v)\le \Delta-1$, then $|\pbar(v)\cap (\pbar(a)\cup \pbar(b))| \le 1$. 
Let $t\in  V(G)\setminus \{a,b\}$ such that $d(t)$ is smallest among the degrees of vertices from $V(G)\setminus \{a,b\}$.  As $|\pbar(a)\cup \pbar(b)|=\Delta-k+2$ and $ |\pbar(t)| \le \Delta-k$,
there exist at least two distinct colors $\alpha,\beta\in (\pbar(a)\cup \pbar(b))\setminus \pbar(t)$.  
By the assumption of this case, we know that every vertex from $V(G)\setminus \{a,b,t\}$ has degree at least $\Delta -k+1$. 
Thus 
by facts (a) and (b), we can find  $x,y\in V(G)\setminus \{a,b,t\}$ such that $\alpha\in \pbar(x)$
and $y\in \pbar(y)$. Clearly, $x$ and $y$ satisfy the property in $(*)$.  This completes the proof. 
\qed 

\section{Proof of Lemma~\ref{lem:kite}}

In this section, we prove Lemma~\ref{lem:kite}. 
We start with some notation. 
Let $G$ be a graph and 
$\varphi\in \CC^k(G-e)$ for some edge $e\in E(G)$ and some integer $k\ge 0$. 
For all the concepts below, when we use them later on, 
if we skip $\varphi$, we mean the concept is defined with respect to the current edge coloring.

Let  $x,y\in V(G)$, and  $\alpha, \beta, \gamma\in [1,k]$ be three colors.   
Let $P$ be an 
$(\alpha,\beta)$-chain of $G$ with respect to $\varphi$ that contains both $x$ and $y$. 
If $P$ is a path, denote by $\mathit{P_{[x,y]}(\alpha,\beta, \varphi)}$  the subchain  of $P$ that has endvertices $x$
and $y$.  By \emph{swapping  colors} along  $P_{[x,y]}(\alpha,\beta,\varphi)$, we mean 
exchanging the two colors $\alpha$
and $\beta$ on the path $P_{[x,y]}(\alpha,\beta,\varphi)$.

Define  $P_x(\alpha,\beta,\varphi)$ to be an $(\alpha,\beta)$-chain or an $(\alpha,\beta)$-subchain of $G$ with respect to $\varphi$ that starts at $x$  and ends at a different vertex missing exactly one of $\alpha$ and $\beta$.    
If $x$ is an endvertex of the $(\alpha,\beta)$-chain that contains $x$, then $P_x(\alpha,\beta,\varphi)$ is unique.  Otherwise, we take one segment of the whole chain to be 
$P_x(\alpha,\beta,\varphi)$. We will specify the segment when it is used.

If  $u$  is a vertrex on  $P_x(\alpha,\beta,\varphi)$, we  write  {$\mathit {u\in P_x(\alpha,\beta, \varphi)}$}; and if $uv$  is an edge on  $P_x(\alpha,\beta,\varphi)$, we  write  {$\mathit {uv\in P_x(\alpha,\beta, \varphi)}$}.  
If $u,v\in P_x(\alpha,\beta,\varphi)$ such that $u$ lies between $x$ and $v$ on the path, 
then we say that $P_x(\alpha,\beta,\varphi)$ \emph{meets $u$ before $v$}. 
Suppose the current color of an  edge $uv$ of $G$
is $\alpha$, the notation  $\mathit{uv: \alpha\rightarrow \beta}$  means to recolor  the edge  $uv$ using the color $\beta$.  Suppose that   $\alpha\in \pbar(x)$ and  $\beta,\gamma\in \varphi(x)$. An $\mathit{(\alpha,\beta)-(\beta,\gamma)}$
\emph{swap at $x$}  consists of two operations:  first swaps colors on $P_x(\alpha,\beta, \varphi)$ to get an edge  $k$-coloring $\varphi'$, and then swaps
colors on $P_x(\beta,\gamma, \varphi')$. 
By convention, an	$(\alpha,\alpha)$-swap at $x$ does  nothing at $x$.

Let $\alpha, \beta, \gamma, \tau,\eta\in[1,k]$. 
We will use a  matrix with two rows to denote a sequence of operations  taken  on $\varphi$.
Each entry in the first row represents a path or  a sequence of vertices. 
Each entry in the second row, indicates the action taken on the object above this entry. 
We require the operations to be taken to follow the ``left to right'' order as they appear in 
the matrix. 
For example,   the matrix below indicates 
three sequential operations taken on the graph based 
on the coloring from the previous step:
\[
\begin{bmatrix}
P_{[a, b]}(\alpha, \beta) &   rs & ab \\
\alpha/\beta & \gamma \rightarrow \tau & \eta
\end{bmatrix}.
\]
\begin{enumerate}[Step 1]
	\item Swap colors on the $(\alpha,\beta)$-subchain $P_{[a, b]}(\alpha, \beta,\varphi) $.
	
	\item Do  $rs: \gamma \rightarrow \tau $. 
	\item Color the edge $ab$ using color $\eta$. 
\end{enumerate}

\begin{figure}[!htb]
	\begin{center}
		\begin{tikzpicture}[scale=1,rotate=90]
			
			{\tikzstyle{every node}=[draw ,circle,fill=white, minimum size=0.5cm,
				inner sep=0pt]
				\draw[blue,thick](0,-2) node (a)  {$a$};
				\draw[blue,thick](0,-3.5) node (b)  {$b$};
				
				\draw [blue,thick](0, -5) node (u)  {$u$};
				\draw [blue,thick](0, -6.5) node (s)  {$s$};
				\draw [blue,thick](0, -8) node (t)  {$t$};
			}
			\path[draw,thick,black!60!green]

			(b) edge node[name=la,pos=0.7, above] {\color{blue} $\alpha$\quad\quad} (u)
			(u) edge node[name=la,pos=0.7, above] {\color{blue}$\beta$\quad\quad} (s)
			(s) edge node[name=la,pos=0.7,above] {\color{blue} $\gamma$\quad\quad} (t);
			

			\draw[dashed, red, line width=0.5mm] (b)--++(140:1cm); 
			\draw[dashed, red, line width=0.5mm] (t)--++(200:1cm); 
			\draw[dashed, red, line width=0.5mm] (t)--++(340:1cm);
			\draw[dashed, red, line width=0.5mm] (a)--++(40:1cm); 
			
			\draw[dashed, red, line width=0.5mm] (a)--++(140:1cm);

			\draw[blue] (-0.5, -3.4 ) node {$\beta$};
			\draw[blue] (0.6, -1.8) node {$\gamma$};
			\draw[blue] (-0.6, -1.8) node {$\alpha$};
			\draw[blue] (0.6, -8.) node {$\beta$};
			\draw[blue] (-0.6, -8.) node {$\alpha$};
			
			{\tikzstyle{every node}=[draw, red ,circle,fill=red, minimum size=0.05cm,
				inner sep=0pt]
				\draw(-0.2,-1.4) node (f1)  {};
				\draw(0,-1.4) node (f1)  {};
				\draw(0.2, -1.4) node (f1)  {};
				
			} 
			
		\end{tikzpicture}

	\end{center}
	\caption{Colors on a Kierstead path of 5 vertices}
	\label{f3}
\end{figure}
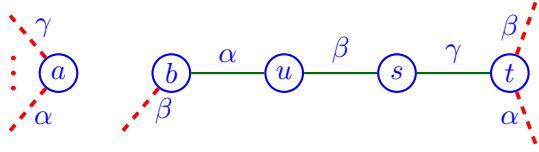

 We will need the following result on Kierstead paths with 5 vertices in proving 
 Lemma~\ref{lem:kite}. 
General properties on Kierstead paths with 5 vertices was proved by the first author 
of this paper~\cite{K5}. Here we stress only one of the cases. 
\begin{LEM}\label{lem:5vexKpathsettingup}
	Let $G$ be a class 2 graph, $ab\in E(G)$ be a critical edge,  $\varphi\in \CC^\Delta(G-ab)$, and  
	$K=(a, ab,b,bu,u, us, s, st, t)$ be a Kierstead path with respect to $ab$ and $\varphi$. 
	If $|\pbar(t)\cap(\pbar(a)\cup \pbar(b))|\ge 3$, then the following holds:
	\begin{enumerate}[(a)]
		\item 
		There exists $\varphi^*\in \CC^\Delta(G-ab)$
		satisfying  the following properties:
		\begin{enumerate}[(i)]
			\item  $\varphi^*(bu)\in \pbar^*(a)\cap \pbar^*(t)$, 
			\item  $\varphi^*(us)\in \pbar^*(b)\cap \pbar^*(t)$, and 
			\item $\varphi^*(st)\in \pbar^*(a)$. 
		\end{enumerate}
		\item  $d_G(b)=d_G(u)=\Delta$. 
		
	\end{enumerate}
	Figure~\ref{f3} shows a Kierstead path with the properties described in (a). 
\end{LEM}

\pf We prove (a). By Lemma~\ref{vf}, for every $\varphi'\in \CC^\Delta(G-ab)$,  $\{a,b\}$ is $\varphi'$-elementary and for every 
$i\in \pbar'(a)$ and $j\in \pbar'(b)$, $a$ and $b$ are $(i,j)$-linked with respect to $\varphi'$. 

Let $\Gamma=\pbar(t)\cap(\pbar(a)\cup \pbar(b))$, and 
$\alpha,\beta \in \Gamma$ be distinct. 
If $\alpha, \beta \in \pbar(a)$, then we let $\lambda\in \pbar(b)$, and do a $(\beta,\lambda)$-swap 
at $b$. If $\alpha, \beta \in \pbar(b)$, then we let $\lambda\in \pbar(a)$, and do a $(\beta,\lambda)$-swap 
at $a$.  Therefore, 
we may assume that $$\alpha\in \pbar(a) \quad \text{and} \quad \beta\in \pbar(b).$$

If $\varphi(bu)=\delta\ne \alpha$, which implies $\delta\in \pbar(a)$, 
then we do a $(\beta,\delta)$-swap and then an $(\alpha,\beta)$-swap at $t$, 
and rename the color $\delta$ as $\alpha$ and vice versa. 
Thus we may assume $$\varphi(bu)=\alpha.$$

Assume first that $\varphi(us)\in \pbar(b)$ and let $\varphi(us)=\tau$. 
If $us\not\in P_t(\beta, \tau)$, we do a $(\beta,\tau)$-swap at $t$
and still call the resulting coloring by $\varphi$, we see that $\tau\in \pbar(b)\cap \pbar(t)$. By exchanging the colors $\beta$ and $\tau$, we have $\varphi(us)=\beta$.
If $us\in P_t(\beta, \tau)$, we do a $(\beta,\tau)$-swap at $t$ and further do an $(\alpha,\beta)$-swap and then an $(\alpha,\tau)$-swap at $t$, which 
gives a new coloring, still call it $\varphi$, such that $\varphi(us)=\beta$.
Let $\varphi(st)=\gamma$.  Since $\alpha,\beta\in \pbar(t)$, $\gamma\ne \alpha, \beta$. 
If $\gamma\in \pbar(a)$, we are done. So we assume  $\gamma\in \pbar(b)\cup \pbar(u)$.  
We color $ab$ by $\alpha$ and uncolor $bu$. Denote this resulting coloring by $\varphi'$.
Then $K'=(b,bu, u, us,s,st,t)$ is a Kierstead path with respect to $bu$ and $\varphi'$. 
However, $\alpha,\beta\in \pbar'(t)\cap (\pbar'(b)\cup \pbar'(u))$, showing a contradiction to 
Lemma~\ref{Lemma:kierstead path1} (b).

Thus we let $\varphi(us)=\delta\in \pbar(a)$. Then $\delta\ne \beta$
by Lemma~\ref{vf} (a).
Let $\varphi(st)=\gamma$. Clearly, $\gamma\ne \alpha,\beta, \delta$. We have either $\gamma\in \pbar(a)$
or $\gamma\in \pbar(b)\cup \pbar(u)$. 
We consider three cases below. 

{\bf \noindent Case 1. $\gamma\in \pbar(b)$.}

If $u\in P_a(\beta,\delta)=P_b(\beta,\delta)$,  we do a $(\beta,\delta)$-swap at $t$.
Since $a$ and $b$ are $(\delta,\gamma)$-linked and $u\in P_t(\delta,\gamma)$,  
we do a $(\delta,\gamma)$-swap at $a$. By renaming colors, this gives a desired coloring $\varphi^*$. 

If $u\not\in P_a(\beta,\delta)=P_b(\beta,\delta)$, we first do a $(\beta,\delta)$-swap at $a$ and then a $(\beta,\gamma)$-swap at $a$. Again this gives a desired coloring $\varphi^*$ after renaming colors.   

\smallskip 

{\bf \noindent Case 2. $\gamma\in \pbar(u)$.}

If $\delta\in \Gamma$, since $b$ and $u$
are $(\beta,\gamma)$-linked by Lemma~\ref{vf} (b),  we do $(\beta,\gamma)$-swap at $t$. 
Now $u\in P_t(\delta,\beta)$, we  do a $(\beta,\delta)$-swap at $a$. 
This gives a desired coloring $\varphi^*$. 
Thus we assume $\delta\not\in \Gamma$. Since $b$ and $u$
are $(\beta,\gamma)$-linked by Lemma~\ref{vf} (b), 
and $a$ and $u$ are $(\delta,\gamma)$-linked by Lemma~\ref{thm:vizing-fan2}~\eqref{thm:vizing-fan2-a}, 
we do $(\beta,\gamma)-(\gamma,\delta)$-swaps at $t$. 
Finally, since $u\in P_t(\beta,\delta)$, we do a $(\beta,\delta)$-swap at $a$.
This gives a desired coloring $\varphi^*$.

\medskip 

{\bf \noindent Case 3. $\gamma\in \pbar(a)$.}

If $\delta\in \Gamma$, we do a $(\beta,\gamma)$-swap at $t$ and then a $(\beta,\delta)$-swap at $a$ to get a desired coloring $\varphi^*$. 
Thus we assume $\delta\not\in \Gamma$. Let $\tau\in \Gamma\setminus \{\alpha,\beta\}$. 
If $\tau\in \pbar(u)$, since $a$ and $u$ are $(\tau, \delta)$-linked by Lemma~\ref{thm:vizing-fan2}~\eqref{thm:vizing-fan2-a}, we do a $(\tau,\delta)$-swap at $t$. This gives back to the previous case that 
$\delta\in \Gamma$. Next we assume $\tau\in \pbar(b)$. 
It is clear that $u\in P_a(\tau,\delta)=P_b(\tau,\delta)$, as otherwise, a $(\tau,\delta)$-swap at $a$ gives a desired coloring. 
Thus we do a $(\tau,\delta)$-swap at $t$,  giving back to the previous case that 
$\delta\in \Gamma$.

Now we assume  $\tau\in \pbar(a)$. 
If $u\not\in P_a(\beta,\delta)$, we do a $(\beta,\delta)$-swap at $a$. Since $a$ and $b$ are $(\alpha,\delta)$-linked and $u\in P_a(\alpha,\delta)$, we do an $(\alpha,\delta)$ swap at $t$.
Now since $u\in P_t(\gamma,\delta)$, we do a  
$(\gamma,\delta)$-swap at $a$, and do  $(\beta,\gamma)-(\gamma,\alpha)$-swaps at $t$.
Since $a$ and $b$ are $(\tau,\gamma)$-linked, we do a 
$(\tau,\gamma)$-swap at $t$, and then a $(\beta,\gamma)$
-swap at $a$. Now since $u\in P_t(\beta, \delta)$, 
we do a $(\beta, \delta)$-swap at $a$. This gives a desired coloring. 
Thus, we assume $u\in P_a(\beta,\delta)$. We do a $(\beta,\delta)$-swap at $t$, and then a $(\tau,\beta)$-swap at $t$.  
Next we do a $(\beta,\gamma)$-swap at $a$ and then a $(\gamma,\delta)$-swap at $a$. 
This gives a desired coloring. 


For statement (b), let $\varphi^*\in \CC^\Delta(G-ab)$ satisfying (i)--(iii). 
Let $\alpha,\gamma\in \pbar^*(a)$, $\beta\in \pbar^*(b)$  with $\alpha,\beta \in \pbar(t)$ such that 
$$
\varphi^*(bu)=\alpha, \quad \varphi^*(us)=\beta,\quad \text{and} \quad  \varphi^*(st)=\gamma. 
$$
Let  $\tau\in \pbar^*(t)\setminus\{\alpha,\beta\}$. 
Suppose to the contrary first that $d_G(b)\le \Delta-1$. Let $\lambda\in \pbar^*(b)\setminus\{\beta\}$. 
We do $(\tau,\lambda)-(\lambda,\gamma)$-swaps at $t$. 
Now we color $ab$ by $\alpha$ and uncolor $bu$ to get a coloring $\varphi'$. 
Then $K'=(b,bu, u, us, s, st, t)$ is a Kierstead path with respect to $bu$ and $\varphi'$. However, $\alpha,\beta \in \pbar'(t)\cap (\pbar'(b)\cup \pbar'(u))$, contradicting Lemma~\ref{Lemma:kierstead path1} (b). 

Assume then that $ d_G(b)=\Delta$ and $d_G(u)\le \Delta-1$. Let $\lambda\in \pbar^*(u)$.
Since $(a,ab,b,bu,u)$ is a multifan, $\lambda\not\in \{\alpha,\beta,\gamma\}$.
Since $u$ and $b$ are $(\beta,\lambda)$-linked and $u$ and $a$ are 
$(\gamma,\lambda)$-linked by Lemma~\ref{thm:vizing-fan2}~\eqref{thm:vizing-fan2-b},  
we do $(\beta,\lambda)-(\lambda,\gamma)$-swap(s) at $t$. 
Now we color $ab$ by $\alpha$ and uncolor $bu$ to get a coloring $\varphi'$. 
Then $K'=(b,bu, u, us, s, st, t)$ is a Kierstead path with respect to $bu$ and $\varphi'$. However, $\alpha \in \pbar'(t)\cap \pbar'(u)$, which contradicts    Lemma~\ref{Lemma:kierstead path1} (a) since $d_G(u)<\Delta$. 
\qed 

We are now ready to prove Lemma~\ref{lem:kite}.
\begin{LEM4}
	Let $G$ be a class 2 graph, $H\subseteq G$ 
	be a kite with $V(H)=\{a,b,c,u,s_1,s_2,t_1,t_2\}$, and let $\varphi\in \CC^\Delta(G-ab)$. 
	Suppose $$K=(a,ab,b,bu,u,us_1, s_1,s_1t_1,t_1) \quad \text{and} \quad K^*=(b,ab,a,ac,c,cu,u,us_2, s_2,s_2t_2,t_2)$$
	are two Kierstead paths with respect to $ab$ and $\varphi$.
	If $\varphi(s_1t_1)=\varphi(s_2t_2)$, 
	then $|\pbar(t_1)\cap \pbar(t_2) \cap ( \pbar(a)\cup \pbar(b))|\le 4$.  
\end{LEM4}

\pf Let $\Gamma=\pbar(t_1)\cap \pbar(t_2) \cap ( \pbar(a)\cup \pbar(b))$. 
Assume to the contrary that $|\Gamma|\ge 5$. By considering $K$ and applying Lemma~\ref{lem:5vexKpathsettingup}, we conclude that $d_G(b)=d_G(u)=\Delta$.
We show that there exists $\varphi^*\in \CC^\Delta(G-ab)$ satisfying the following properties:
\begin{enumerate}[(i)]
	\item  $\varphi^*(bu), \varphi^*(cu), \varphi^*(us_2)\in \pbar^*(a)\cap \pbar^*(t_1)\cap \varphi^*(t_2)$, 
	\item  $\varphi^*(us_1)\in \pbar^*(b)\cap \pbar^*(t_1)\cap \varphi^*(t_2)$, and 
	\item $\varphi^*(s_1t_1)=\varphi^*(s_2t_2)\in \pbar^*(a)$. 
\end{enumerate}

\begin{figure}[!htb]
	\begin{center}
		\begin{tikzpicture}[scale=1]
			
			{\tikzstyle{every node}=[draw ,circle,fill=white, minimum size=0.5cm,
				inner sep=0pt]
				\draw[blue,thick](0,-2) node (a)  {$a$};
				\draw[blue,thick](-1,-3) node (b)  {$b$};
				\draw[blue,thick](1,-3) node (c)  {$c$};
				\draw [blue,thick](0, -4.5) node (u)  {$u$};
				\draw [blue,thick](-1, -6) node (x)  {$s_1$};
				\draw [blue,thick](1, -6) node (y)  {$s_2$};
				\draw [blue,thick](-1, -7.5) node (t1)  {$t_1$};
				\draw [blue,thick](1, -7.5) node (t2)  {$t_2$};
			}
			\path[draw,thick,black!60!green]
			(a) edge node[name=la,pos=0.8, above] {\color{blue} $\beta$} (c)
			
			(c) edge node[name=la,pos=0.4, below] {\color{blue} \quad$\tau$} (u)
			(b) edge node[name=la,pos=0.4, below] {\color{blue} $\alpha$\quad\quad} (u)
			(u) edge node[name=la,pos=0.6, above] {\color{blue}$\beta$\quad\quad} (x)
			(x) edge node[name=la,pos=0.7, above] {\color{blue}$\gamma$\quad\quad} (t1)
			(y) edge node[name=la,pos=0.7,above] {\color{blue}  $\gamma$\quad\quad} (t2)
			(u) edge node[name=la,pos=0.6,above] {\color{blue}  \quad$\delta$} (y);
			

			\draw[dashed, red, line width=0.5mm] (b)--++(140:1cm); 
			\draw[dashed, red, line width=0.5mm] (t1)--++(200:1cm); 
			\draw[dashed, red, line width=0.5mm] (t1)--++(250:1cm); 
			\draw[dashed, red, line width=0.5mm] (t1)--++(290:1cm); 
			\draw[dashed, red, line width=0.5mm] (t1)--++(340:1cm); 
			\draw[dashed, red, line width=0.5mm] (t2)--++(200:1cm); 
			\draw[dashed, red, line width=0.5mm] (t2)--++(250:1cm); 
			\draw[dashed, red, line width=0.5mm] (t2)--++(290:1cm); 
			\draw[dashed, red, line width=0.5mm] (t2)--++(340:1cm); 
			
			\draw[dashed, red, line width=0.5mm] (a)--++(40:1cm); 
			\draw[dashed, red, line width=0.5mm] (a)--++(100:1cm); 
			\draw[dashed, red, line width=0.5mm] (a)--++(70:1cm); 
			\draw[dashed, red, line width=0.5mm] (a)--++(140:1cm); 

			\draw[blue] (-1.6, -9+1.5) node {$\alpha$}; 
			\draw[blue] (1.6, -9+1.5) node {$\alpha$}; 
			\draw[blue] (-1.4, -9.5+1.5) node {$\beta$}; 
			\draw[blue] (1.4, -9.5+1.5) node {$\beta$}; 
			\draw[blue] (-1, -9.6+1.5) node {$\tau$}; 
			\draw[blue] (1, -9.6+1.5) node {$\tau$}; 
			\draw[blue] (-0.45, -9.4+1.5) node {$\delta$}; 
			\draw[blue] (0.45, -9.4+1.5) node {$\delta$};

			\draw[blue] (-1.2, -2.5) node {$\beta$};
			\draw[blue] (0.6, -1.8) node {$\gamma$};
			\draw[blue] (-0.6, -1.8) node {$\alpha$};
			\draw[blue] (-0.3, -1.4) node {$\tau$};
			\draw[blue] (0.4, -1.4) node {$\delta$};
			
			%
			%
			
		\end{tikzpicture}
	\end{center}
	\caption{Colors on the edges of a kite}
	\label{pic2}
\end{figure}
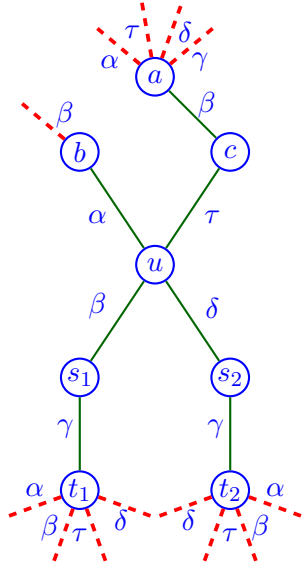

See Figure~\ref{pic2} for a depiction of the colors described above. 

Let $\alpha,\beta, \tau,\delta\in \Gamma$ be distinct, and let $\varphi(s_1t_1)=\varphi(s_2t_2)=\gamma$. 
We may assume that $\alpha\in \pbar(a)$ and $\beta\in \pbar(b)$. 
Otherwise, since $d_G(b)=\Delta$, we have $\alpha,\beta \in \pbar(a)$. 
Let $\lambda\in \pbar(b)$. As $a$ and $b$ are  $(\beta,\lambda)$-linked, we do a  
$(\beta,\lambda)$-swap at $b$. Note that this operation may change some colors of the edges of $K$
and $K^*$, but they are still Kierstead paths with respect to $ab$
and the current coloring.  Note that $\varphi(ac)=\beta$, since $\pbar(b)=\{\beta\}$ and $K^*$ is a Kierstead path.

Since $d_G(b)=d_G(u)=\Delta$, and $\beta\in \pbar(b)\cap \pbar(t_1)$, 
we know that $\gamma\in \pbar(a)$, as  $K_1$ is a
Kierstead path. 
Next, we may assume that $\varphi(bu)=\alpha$. 
If not, let $\varphi(bu)=\alpha'$. Since $a$
and $b$ are $(\alpha,\beta)$-linked, we do an $(\alpha,\beta)$-swap at $b$.
Now $a$ and $b$ are   $(\alpha,\alpha')$-linked, we do an  $(\alpha,\alpha')$-swap at $b$. 
Finally, we do an $(\alpha',\beta)$-swap at $b$. 
All these swaps do not change the colors in $\Gamma$, so now we get the color on $bu$
to be $\alpha$. 

We may now assume that $\varphi(cu)=\tau$. 
If not, let $\varphi(cu)=\tau'$.  Then $\tau'\in \pbar(a)$. 
Since $a$ and $b$ are $(\beta,\tau)$-linked, we do a $(\beta,\tau)$-swap at $b$. 
Then do  $(\tau,\tau')-(\tau',\beta)$-swaps at $b$. 

Finally, we show that we can modify $\varphi$
to get $\varphi'$ such that $\varphi'(us_1)=\beta$
and $\varphi'(us_2)=\delta$. Let $\lambda\in \Gamma\setminus \{\alpha,\beta,\tau,\delta\}$.
Assume firstly that $\varphi(us_1)=\beta'\ne \beta$. Then $\beta'\in \pbar(a)$ as $\pbar(b)=\{\beta\}$ and $\pbar(u)=\emptyset$. 
If $\beta'\in \Gamma$, we do $(\beta,\gamma)-(\gamma,\beta')$-swaps at $b$. 
Thus, we assume $\beta'\not\in \Gamma$. 
If $u\not\in P_a(\beta,\beta')=P_b(\beta,\beta')$, we simply do a $(\beta,\beta')$-swap at $b$. 
Thus, we assume $u\in P_a(\beta,\beta')=P_b(\beta,\beta')$. 
We do a $(\beta,\beta')$-swap at both $t_1$
and $t_2$. Since $a$
and $b$ are $(\beta,\lambda)$-linked, we 
do a $(\beta,\lambda)$-swap at both $t_1$
and $t_2$. 
Now we do $(\beta,\gamma)-(\gamma,\beta')$-swaps at $b$. 
By switching the role of $\beta$ and $\beta'$, 
we have $\varphi(us_1)=\beta$.

Lastly, we show that $\varphi(us_2)=\delta$. 
Assume otherwise that $\varphi(us_2)=\delta'$. 
By coloring $ab$ with the color $\beta$ and uncoloring $ac$, 
we have $d_G(c)=\Delta$ by Lemma~\ref{lem:5vexKpathsettingup}.
Thus $\delta'\in \pbar(a)$.  
Note that $bu\in P_{t_1}(\alpha,\gamma)$. 
Otherwise, let $\varphi'=\varphi/P_{t_1}(\alpha,\gamma)$.
Then $P_b(\alpha,\beta)=bus_1t_1$, showing 
a contradiction to the fact that $a$
and $b$ are $(\alpha,\beta)$-linked with respect to $\varphi'$. 
Thus, $bu\in P_{t_1}(\alpha,\gamma)$.
Next, we claim that $P_{t_1}(\alpha,\gamma)$ meets $u$ before $b$. 
As otherwise, we do the following operations to get a $\Delta$-coloring of $G$:
\[
\begin{bmatrix}
	s_1t_1& P_{[s_1,b]}(\alpha,\gamma)  &  us_1 &  bu  & ab\\
	\gamma\rightarrow \beta& \alpha/\gamma &  \beta \rightarrow \alpha&
	\alpha\rightarrow \beta&  \gamma\end{bmatrix}.
\]
This gives a contradiction to the assumption that $G$ is $\Delta$-critical. 
Thus, we  have that $P_{t_1}(\alpha,\gamma)$ meets $u$ before $b$.
This implies that it is not the case that $P_{t_2}(\alpha,\gamma)$
meets $u$ before $b$. In turn, this implies that $u\in P_a(\beta,\delta')=P_b(\beta,\delta')$.  
As otherwise,  we get a $\Delta$-coloring of $G$ by doing a $(\beta,\delta')$-swap along the $(\beta,\delta')$-chain containing $u$,
and then doing the same operation as above with $t_2$
playing the role of $t_1$. 

Since $u\in P_a(\beta,\delta')=P_b(\beta,\delta')$, 
we do a $(\beta,\delta')$-swap at both $t_1$ and $t_2$. As $u\in P_a(\beta,\tau)=P_b(\beta,\tau)$, 
we do a $(\beta,\tau)$-swap at both $t_1$
and $t_2$. Since $us_1\in P_{t_1}(\beta,\gamma)$, 
we do a $(\beta,\gamma)$-swap at $b$, then a  $(\gamma,\lambda)$-swap at $b$. 
Since $a$ and $b$ are $(\tau,\lambda)$-linked, we do a $(\tau,\lambda)$-swap at both $t_1$
and $t_2$. Now $(\lambda,\delta)-(\delta,\gamma)-(\gamma,\beta)$-swaps at $b$
give a desired coloring. 

Still, by the same arguments as above,  we have that  $P_{t_1}(\alpha,\gamma)$ meets $u$ before $b$,
and   $u\in P_a(\beta,\delta)=P_b(\beta,\delta)$. 
Let $P_u(\beta,\delta)$ be the $(\beta,\delta)$-chain starting at $u$ not including the edge $us_2$. 
It is clear that $P_u(\beta,\delta)$  ends at either $a$ or $b$. 
We may assume that $P_u(\beta,\delta)$ ends at $a$. 
Otherwise, we color $ab$ by $\beta$, uncolor $ac$, and let $\tau$
play the role of $\alpha$. Let $P_u(\alpha,\gamma)$
be the $(\alpha,\gamma)$-chain starting at $u$ not including the edge $bu$,
which ends at $t_1$ by our earlier argument. 
We do the following operations to get a $\Delta$-coloring of $G$:
\[
\begin{bmatrix}
	P_u(\alpha,\gamma)& bu &  P_u(\beta,\delta) &  us_2t_2 & ab\\
	\alpha/\gamma& \alpha\rightarrow \beta &  \beta/\delta&
	\delta/\gamma&  \alpha\end{bmatrix}.
\]
This gives a contradiction to the assumption that $G$ is $\Delta$-critical.
The proof is now finished. 
\qed


\end{document}